\documentclass[12pt]{amsart}
\usepackage{amssymb}
\usepackage{amscd}
\usepackage{pstricks}
\usepackage{pst-node}
\newtheorem{thm}{Theorem}[section]
\newtheorem{cor}[thm]{Corollary}

\newtheorem{prop}[thm]{Proposition}

\theoremstyle{definition}

\newtheorem{rem}[thm]{Remark}

\theoremstyle{remark}

\numberwithin{equation}{section}

\newcommand{\GL}{\mathit{GL}}

\newcommand{\trans}[1]{{}^t\kern-.2em{#1}}
\newcommand{\ytrans}[1]{{}^t\kern-.11em{#1}}
\newcommand{\Trans}[1]{{}^T\kern-.2em{#1}}
\newcommand{\lsup}[2]{{}^{#1}\kern-.1em{#2}}

\newcommand{\Ker}{\operatorname{Ker}}
\newcommand{\Coker}{\operatorname{Coker}}

\newcommand{\M}{\mathbf{\M}}

\DeclareFixedFont{\bgn}{OT1}{cmr}{m}{n}{20.74}
\DeclareFixedFont{\bgi}{OT1}{cmr}{m}{it}{20.74}

\newcommand{\bigzerou}{\smash{\lower1.7ex\hbox{\bgi O}}}

\makeatletter
\def\eqnarray{%
   \stepcounter{equation}%
   \def\@currentlabel{\p@equation\theequation}%
   \global\@eqnswtrue
   \m@th
   \global\@eqcnt\z@
   \tabskip\@centering
   \let\\\@eqncr
   $$\everycr{}\halign to\displaywidth\bgroup
       \hskip\@centering$\displaystyle\tabskip\z@skip{##}$\@eqnsel
      &\global\@eqcnt\@ne \hfil$\displaystyle{{}##{}}$\hfil
      &\global\@eqcnt\tw@ $\displaystyle{##}$\hfil\tabskip\@centering
      &\global\@eqcnt\thr@@ \hb@xt@\z@\bgroup\hss##\egroup
        \tabskip\z@skip
      \cr}
\makeatother
\makeatletter
\def\varin{\mathrel{\mathpalette\@varin\relax}}
\def\@varin#1{%
   \hbox{\setbox\z@\hbox{\m@th$#1\cup$}%
       \def\reserved@a{bold}%
       \dimen@\ifx\reserved@a\math@version .3\else .2\fi\p@
       \kern.5\wd\z@\kern-\dimen@
       \vrule\@width2\dimen@\@height1.08\ht\z@\@depth\z@
       \kern-\dimen@\kern-.5\wd\z@
       \box\z@}}
\makeatother

\begin{document}
 


\title [Quillen determinant] 
{On the Quillen determinant}



\author{Kenro Furutani}
\address{Kenro Furutani \endgraf 
Department of Mathematics \endgraf 
Faculty of Science and Technology \endgraf 
Science University of Tokyo \endgraf 
2641 Noda, Chiba (278-8510)\endgraf 
Japan \endgraf}
\email{furutani@ma.noda.sut.ac.jp}

\bigskip
\bigskip

\begin{abstract}
We explain the bundle structures of the {\it Determinant line bundle} 
and the {\it Quillen determinant line bundle} 
considered on the connected component of the space of
Fredholm operators including the identity operator in an intrinsic way.
Then we show that these two are isomorphic and that they are 
non-trivial line bundles and trivial
on some subspaces. Also we remark 
a relation of the {\it Quillen determinant line bundle}  
and the {\it Maslov line bundle}.
\end{abstract}

\maketitle
\tableofcontents
\thispagestyle{empty}

\section{Introduction}

The Fredholm determinant is defined for the class of the operators of the form
``{\it Id + trace class operator}'' on a Hilbert space $H$
as the extension of the finite dimensional cases with respect to
the trace norm :
\begin{equation}
{\det}_{F}(Id + K)= \prod (1+\lambda_i),  
\end{equation}
where $\lambda_i$ are eigenvalues of the trace class operator 
$K$ (see \cite{GGK} for analytic properties of the Fredholm determinant).
This quantity gives us a $\mathbb{C}^*$-valued
holomorphic one cocycle on the space of Fredholm operators on $H$
whose Fredholm indexes are zero. 
In fact,
let $\mathcal{F}=\mathcal{F}(H)$ be the space of Fredholm operators 
defined on a Hilbert space $H$ and 
we denote by $\mathcal{F}_0=\mathcal{F}_0(H)$ 
the connected component of $\mathcal{F}(H)$
consisting of the operators with the index zero. 
Let $\mathcal{I}_1$ be the space of trace class operators on $H$.
For each trace class operator $A\in\mathcal{I}_1$, 
we denote by $\mathcal{U}_A$ an open subset
of $\mathcal{F}_0$ 
consisting of such operators $T$ that
$T+A$ is an isomorphism of $H$. Then $\mathcal{F}_0$ is covered
by open subsets $\{\mathcal{U}_A\}_{A\in\mathcal{I}_1}$.
Let $A$ and $B$ be two trace class operators, and let
$T \in\mathcal{U}_A\cap\mathcal{U}_B \subset \mathcal{F}_0$,
then
the functions $\left\{g_{A,B}(T)\right\}_{A,\,B\in\mathcal{I}_1}$
$$
g_{A,B}(T) ={\det}_F\left\{Id + (A-B)(T+B)^{-1}\right\}
= {\det}_F\left\{(T+A)(T+B)^{-1}\right\}
$$
are holomorphic on $\mathcal{U}_A\cap\mathcal{U}_B$, 
and they satisfy the cocycle condition :
\begin{equation}\label{eq:transition}
g_{A,C}(T)=g_{A,B}(T)g_{B,C}(T)
\end{equation}
for $T\in \mathcal{U}_A\cap\mathcal{U}_B\cap\mathcal{U}_C$.
We denote by $\mathcal{L}_F$ the complex line bundle defined by
these transition functions 
$\{g_{A,B}\}_{A,\,B\in\mathcal{I}_1}$ 
and call it as the ``{\it Determinant line bundle}''.  

The disjoint unions of finite dimensional vector spaces
\begin{center}
$\coprod\limits_{T\in\mathcal{F}_0}\Ker(T)$ \quad and\quad
$\coprod\limits_{T\in \mathcal{F}_0}\Coker(T)$
\end{center}
do not have vector bundle structures. When we consider them on a
compact subset $X$ in $\mathcal{F}_0$, it can be seen that 
the formal difference of these two is 
an element of the $K$-group $K(X)$ by {\it approximating} 
each of these two with suitable
vector bundles which are constructed by a standard method.

On the other hand,
the disjoint union of the lines
$$
\coprod\limits_{T\in\mathcal{F}_0}
\stackrel{\dim \Ker(T)}\wedge\Ker(T)^*\otimes
\stackrel{\dim \Coker(T)}\wedge
\Coker(T)
$$
has a complex line bundle structure on the whole space
$\mathcal{F}_0$ and is called as the {\it Quillen determinant line bundle}.
This fact is stated in the paper \cite{Qu} and it is treated in various
contexts(\cite{At2}, \cite{Ma}, \cite{Pa}, \cite{SW} and others).

In this note we give a rigorous proof of this fact by giving 
an intrinsic correspondence between the Determinant line bundle 
and the Quillen determinant($\S 2$),
and prove that it is non-trivial on $\mathcal{F}_0$($\S 3$).
In $\S 4$ we show it is trivial on each compact subset of the subspace
$\widehat{\mathcal{F}}_*$ (=
the non-trivial connected component in the space of selfadjoint
Fredholm operators). Of course it is trivial on the each 
subspace of essentially positive and essentially negative
Fredholm operators (= $\widehat{\mathcal{F}}_{\pm}$).  
Finally in $\S 5$ we prove that the induced bundle
of the Quillen determinant line bundle on the space
of {\it Fredholm Lagrangian Grassmannian} by a naturally defined map 
is trivial and remark a relation with the {\it Maslov line bundle}.


\section{Fredholm determinant and the Quillen Determinant} 

Let $T$ be a Fredholm operator on a (complex)Hilbert 
space $H$. We denote by $\mathcal{A}_T$ a subset
of the space of trace class operators $\mathcal{I}_1$ such that
$$
\mathcal{A}_T=\{A\in \mathcal{I}_1~|~ T+A ~\text{is invertible}\},
$$
and let $\mathcal{D}_T$ be a space of complex valued functions 
on $\mathcal{A}_T$
satisfying the following condition :
$$
\mathcal{D}_T=
\left\{f:\mathcal{A}_T\to \mathbb{C}\,|\, 
f(B)={\det}_F\{(T+A)(T+B)^{-1}\}f(A)\right\}.
$$

Of course this is an 1-dimensional vector space, and the union 
$$
\coprod\limits_{T\in\mathcal{F}_0(H)}\mathcal{D}_T
$$
becomes a holomorphic complex line bundle with local trivializations
\begin{align}
j_A\,:\,&\coprod\limits_{T\in\mathcal{U}_A}\mathcal{D}_T
\stackrel{\sim}\to \mathcal{U}_A\times \mathbb{C}\label{l-t-1}\\
j_A\,:\,\mathcal{D}_T\ni f
&\mapsto (T,f(A))\in \mathcal{U}_A\times\mathbb{C},\notag
\end{align}
where $A\in\mathcal{I}_1$ and 
$\mathcal{U}_A=\{ T\in\mathcal{F}_0(H)\,|\,~T+A ~\text{is invertible}\}$.
By the definition of the function space $\mathcal{D}_T$,
the transition function on $\mathcal{U}_A\bigcap\mathcal{U}_B$
is given by
$$
{\det}_F\left\{(T+A)(T+B)^{-1}\right\},
$$
so that the space
$$
\coprod\limits_{T\in\mathcal{F}_0(H)}\mathcal{D}_T
$$
is a realization of the {\it Determinant line bundle}.
We denote it by $\mathcal{L}_F$.

For a fixed $T\in \mathcal{F}_0$, 
we denote by $\pi_T :H\to H$ the orthogonal
projection operator onto the $\Ker(T)$ and by $\rho_T$ the natural
projection $\rho_T : H \to \Coker(T)$. 

Let $L$ be a linear map $L : Ker(T) \to H$ satisfying the condition :
\begin{equation}\label{condition-L}
The ~composition ~\rho_T\circ L :\Ker(T) \to \Coker(T)
~is ~an ~isomorphism. 
\end{equation}
Then under this condition for the operator $L$ we know that 
the operator $T+L\circ\pi_T$ is an isomorphism on $H$. 

Let $\{e_i\}_{i=1}^{d}$ 
be a basis of $\Ker(T)$ and $\{e_i^{\,*}\}_{i=1}^d$ 
the dual basis ($d=\dim\Ker T$). 
Then we define a map
\begin{equation}\label{def:natural-map}
\phi_T: \mathcal{D}_T\, \to \stackrel{\dim \Ker(T)}\wedge\Ker(T)^*
\otimes
\stackrel{\dim\Coker(T)}\wedge
\Coker(T)
\end{equation}
by
\begin{align}
&\phi_T(f)\label{l-t-2}\\ 
&=f(A)\cdot{\det}_F\left\{(T+A)(T+L\circ\pi_T)^{-1}\right\}\times\notag
\\
&\qquad\qquad\qquad\times 
e_1^*\wedge\cdots\wedge e_{d}^*\otimes\rho_T(L(e_1))
\wedge\cdots\wedge\rho_T(L(e_{d})),\notag
\end{align}
where we fixed an $A\in\mathcal{A}_T$. By the relation
\begin{align*}
&f(A){\det}_F\left\{(T+A)(T+L\circ\pi_T)^{-1})\right\}
{\det}_F\left\{((T+B)(T+L\circ\pi_T)^{-1})^{-1}\right\}\\
&=f(A){\det}_F\left\{(T+A)(T+B)^{-1}\right\}=f(B),\\
\end{align*}
it will be clear of the independence of the definition of this map
from the choice of $A\in\mathcal{A}_T$ and the map $\phi_T$ 
is an isomorphism.  
Moreover we have
\begin{prop}
The definition of the map $\phi_T$  
depends neither on the choice of the 
map $L$ satisfying the condition (\ref{condition-L}) above,
nor on the choice of the basis
$\{e_i\}$ of $\Ker(T)$.
\end{prop}
\begin{proof}
Again it would be clear of the independence from the choice
of a basis of $\Ker(T)$. So we only prove the independence
from the choice of the operator $L$. 

Let $L'$ be another such operator
$L': \Ker(T) \to H$ that 
$ \rho_T\circ L': \Ker(T) \to \Coker(T)$ is isomorphic,
then we have
$$
\rho_T(L'(e_j)) = \sum_{i}\,a_{i\,j}\rho_T(L(e_i)).
$$
and 
$$
T+L\circ\pi_T = T+L'\circ\pi_T ~\text{on}~ \Ker(T)^{\perp}. 
$$
Hence
$$
(T+L\circ\pi_T)^{-1}\circ(T+L'\circ\pi_T) - Id 
$$
is a finite rank operator, 
and moreover we have 
$$
{\det}_F\left\{(T+L\circ\pi_T)^{-1}\circ(T+L'\circ\pi_T)\right\}
=\det(a_{i\,j}).
$$
This relation gives us
\begin{align}
&f(A)\cdot{\det}_F\left\{(T+A)(T+L\circ\pi_T)^{-1}\right\}\times\\
&\qquad\qquad\qquad\qquad\times e_1^*\wedge\cdots\wedge e_{d}^*\otimes
\rho_T(L(e_1))\wedge\cdots\wedge\rho_T(L(e_{d}))\notag\\
=&
f(A)\cdot{\det}_F\left\{(T+A)(T+L'\circ\pi_T)^{-1}\right\}
\times\notag\\
&\qquad\qquad\qquad\qquad\times e_1^*\wedge\cdots\wedge e_{d}^*\otimes
\rho_T(L'(e_1))\wedge\cdots\wedge\rho_T(L'(e_{d})),\notag
\end{align}
which proves the independence of the definition of the map $\phi_T$
from the choice of the linear map $L$.
\end{proof}

By this proposition we can introduce (the topology and) 
the local trivialization of the space
$$\coprod\limits_{T\in \mathcal{U}_A}\stackrel{\dim \Ker(T)}\wedge\Ker(T)^*
\otimes\stackrel{\dim\Coker(T)}\wedge
\Coker(T)
$$ 
through the local trivialization (\ref{l-t-1})
and the map $\phi_T$ :
\begin{equation}\label{l-t-3}
\left(\coprod\limits_{T\in\mathcal{U}_A}\phi_T\right) \circ j_A^{\,-1}
\,:\,\mathcal{U}_A \times \mathbb{C}\stackrel{\sim}\to 
\coprod\limits_{T\in \mathcal{U}_A}\stackrel{\dim\Ker(T)}\wedge
\Ker(T)^*
\otimes\stackrel{\dim\Coker(T)}\wedge
\Coker(T).
\end{equation}
Then 
$$\coprod\limits_{T\in\mathcal{F}_0}
\stackrel{\dim\Ker(T)}\wedge
\Ker(T)^*
\otimes\stackrel{\dim\Coker(T)}\wedge
\Coker(T)
$$
becomes
a complex line bundle which is isomorphic to the determinant line
bundle
$\mathcal{L}_F$.  This is 
the ``{\it Quillen determinant line bundle}'' and we denote it by
$\mathcal{L}_Q$.


\section{Non triviality of the Quillen determinant}                   
               
\begin{thm}
The bundle $\mathcal{L}_Q$ is not trivial on the whole 
space $\mathcal{F}_0$.
\end{thm}
\begin{proof}
For a compact Hausdorff space $X$ we know by the famous theorem
\cite{At1} that the reduced $K$-group $\widetilde{K}(X)$ is isomorphic to
the space of homotopy classes $[X, \mathcal{F}_0]$ 
of continuous maps $f: X \to \mathcal{F}_0$ and the correspondence
is given by constructing two vector bundles $E$ and $F$ on $X$
which satisfy the following exact sequence at each point $x\in X$ :
\begin{equation}
0\to \Ker(f(x))\to E_x\to F_x \to \Coker(f(x))\to 0.
\end{equation}
The homotopy class of the map $f$ corresponds the the element
$[E] - [F] \in K(X)$.

Hence we have 
$f^*(\mathcal{L}_Q) 
\cong \stackrel{\dim E}\wedge E^*\otimes\stackrel{\dim F}\wedge F,
$
and so for any line bundle $\ell$ on a compact space $X$
the element 
$[\ell]-[\varepsilon^1]\in\widetilde{K}(X)$ 
($\varepsilon^1 :~\text{one dimentional trivial line bundle}$) 
corresponds to a continuous map
$g :X \rightarrow \mathcal{F}_0$, we have
$\ell^*\cong g^*(\mathcal{L}_Q)$.  Hence we know by taking a suitable
compact space $X$ with $H^2(X,\mathbb{Z})\not= \{0\}$ that $\mathcal{L}_Q$
can not be trivial on the whole space $\mathcal{F}_0(H)$.
\end{proof}

\section{A triviality of the Quillen determinant} 

Although we have proved that the Quillen determinant line bundle
is not trivial on the whole space $\mathcal{F}_0$, it might be 
trivial on a subspace in $\mathcal{F}_0(H)$. For example it is trivial
on the space of essentially positive(negative) Fredholm operators
(=$\widehat{\mathcal{F}}_{\pm}$).

Now let $\widehat{\mathcal{F}}_*$ be the non-trivial connected component
of the selfadjoint Fredholm operators.  Then we have
\begin{thm}
On each compact subset in the space $\widehat{\mathcal{F}}_*$
the  Quillen determinant $\mathcal{L}_Q$ is trivial.
\end{thm}

\begin{proof}
Let $X$ be a compact Hausdorff space and let $f$ be a continuous map,
$f:X\to \widehat{\mathcal{F}}_*$.
It is enough to show that $f^*(\mathcal{L}_Q)$ is trivial. Let
$\Omega\mathcal{F}_0$ be the path space consisting
of paths connecting $Id$ and $-Id$.
Let $\alpha : \widehat{\mathcal{F}}_* \to \Omega \mathcal{F}_0$
be a continuous map given by
\begin{equation}
\alpha(A)(t)= \cos(\pi\,t) +\sin(\pi\, t)\cdot A \in
\Omega\mathcal{F}_0,~t\in [0,1].
\end{equation}
This is a homotopy equivalence (see \cite{AS}).

Now let $S(X)$ be the suspension of $X$, then we 
have a continuous map $h_f: S(X) \to \mathcal{F}_0$ defined by 
$$
h_f(t,x)= \alpha(f (x))(t).
$$
Let $\mathcal{C}: X \to S(X)$ be a map defined as $\mathcal{C}(x) = (1/2,x)\in S(X)$,
then by the definition of the suspension we have
$$
h_f\circ\mathcal{C} = ~{\bf i}\circ f,
$$
where ${\bf i}$ is the inclusion map 
${\bf i} : \widehat{\mathcal{F}}_* \hookrightarrow \mathcal{F}_0$.
Since $\widetilde{K}(S(X)) 
= ~\text{ind-lim}_{n \to\infty}[X, \GL(n,\mathbb{C})]$,
we know that the induced map 
$\mathcal{C}^* :\widetilde{K}(S(X))\to \widetilde{K}(X)$ is trivial.
Hence the induced line bundle
$(h_f\circ\mathcal{C})^*(\mathcal{L}_Q)=f^*(\mathcal{L}_Q)$ 
must be trivial.
\end{proof}

\section{Quillen determinant on the Fredholm Lagrangian Grassmannian} 

In this section we show that the Quillen determinant is trivial, when
it is
pull-backed on the {\it Fredholm Lagrangian Grassmannian} through
an embedding.

First we describe the Fredholm Lagrangian Grassmannian. So,
let $H$ here be a real symplectic Hilbert space. The symplectic form
is non-degenerate in such a sense that
$\omega: H \times H \to \mathbb{R}$ defines the continuous 
isomorphism $\omega^{\#}$ :
\begin{equation}
\omega^{\#} : H \to H^*,~\omega^{\#}(x)(y)=\omega(x,y).
\end{equation}

We do not change the symplectic form $\omega$ once it has been
introduced on a real Hilbert space $H$, but rather freely we can
replace the inner product with a new one whose defining norm is equivalent
to that defined by the original inner product. Especially 
we can assume from the beginning
that the symplectic form $\omega$ is expressed in the form
$\omega(x,y) = <\,J(x),\, y\,>$, where $J$ is an almost complex structure
with the property that $<\,J(x),\,J(y)\,> = < x , y >$, $\,^tJ =-J$, $\,^tJ$ is 
the transpose operator with respect to the (Euclidean)inner 
product $<\cdot,\cdot>$.

Let $\lambda$ be a Lagrangian subspace :
\begin{equation}
\lambda
=\lambda^{\circ}=\{x\in H \,|\, \omega(x,y)=0 ~\text{for any} ~y\in \lambda\},
\end{equation}
and denote by $\mathcal{F}\Lambda_{\lambda}(H)$ the space
of such Lagrangian subspaces $\mu$ that the pair $(\mu,\lambda)$
is a Fredholm pair(see \cite{Ka} for a general theory of Fredholm
pairs and \cite{FO} for particular properties of Fredholm pairs of Lagrangian
subspaces), that is,
\begin{align*}
&(\text{i})\qquad\dim(\lambda\cap\mu)< +\infty\\
&(\text{ii})\quad\lambda+\mu 
~\text{is a closed and finite codimensional subspace in} ~H.
\end{align*}
We call this space as the ``{\it Fredholm Lagrangian Grassmannian}''.
The topology is naturally defined by embedding it into the space of
bounded operators $\mathcal{B}(H)$ on $H$
by the map 
$\mathcal{P}: \mathcal{F}\Lambda_{\lambda}(H)\to \mathcal{B}(H)$,
$\mathcal{P}(\mu)$ is the orthogonal projection operator onto
the space $\mu$ and $\mathcal{F}\Lambda_{\lambda}(H)$ 
becomes an infinite dimensional smooth manifold.
It is known that the fundamental group
$\pi_1(\mathcal{F}\Lambda_{\lambda}(H))$ is $\mathbb{Z}$ and the
isomorphism is given by so called the {\it Maslov index} for each loop.

When we regard the real Hilbert space $H$ as a complex Hilbert space by means of 
the almost complex structure $J$ with the
Hermitian inner product 
$$
(x,y)=<x,y>-\sqrt{-1}<J(x),y>,
$$
we denote it by $H_J$.
 
Each
Lagrangian subspace $\lambda$ defines a {\it real structure} on $H_J$ :
\begin{align*}
\lambda\otimes\mathbb{C}&\stackrel{\sim}\to H_J,\\
 x\otimes 1+y\otimes \sqrt{-1}&\mapsto x+J(y).
\end{align*}

We denote by $\tau_{\lambda}$ the complex
conjugation with respect to a real structure given 
by a Lagrangian subspace $\lambda$ :
\begin{equation}
\tau_{\lambda}(x+J(y)) = x-J(y).
\end{equation}
This is an anti-linear involution on $H=H_J$ and 
$2\mathcal{P}(\lambda)-Id=\tau_{\lambda}$.

Let $\mu\in\mathcal{F}\Lambda_{\lambda}(H)$, then the operator
$$
-\tau_{\mu}\circ\tau_{\lambda}
$$
is a unitary operator $\in \mathcal{U}(H_J)$ with the property
that
$$
Id -\tau_{\mu}\circ\tau_{\lambda}
$$
is a Fredholm operator. We denote the correspondence 
$\mu \mapsto -\tau_{\mu}\circ\tau_{\lambda}$
by 
\begin{equation}
\mathcal{S}_{\lambda}:
\mathcal{F}\Lambda_{\lambda}(H)\to \mathcal{U}_F(H_J),
\end{equation}
where $\mathcal{U}_F(H_J)$ is a space of unitary operators
$U$ on $H_J$ such that $U+Id$ is a Fredholm operator. 

We call the map
$\mathcal{S}_{\lambda}$ the {\it Souriau map} 
(\cite{Le}, \cite{So} and \cite{FO}) which satisfies 
$
\mathcal{S}_{\lambda}(\mu)^*=\mathcal{S}_{\mu}(\lambda).
$
We know that through this map the fundamental groups of
the Fredholm Lagrangian Grassmannian and the space $\mathcal{U}_F(H_J)$
are isomorphic.

Let us denote by $q_{\lambda}$ the map
\begin{align*}
&q_{\lambda} : \mathcal{F}\Lambda_{\lambda}(H) \to \mathcal{F}_0(H_J)\\
&q_{\lambda}(\mu)=Id-\tau_{\mu}\circ\tau_{\lambda}
=Id+\mathcal{S}_{\lambda}(\mu).
\end{align*}

\begin{thm}
The pull back $q_{\lambda}\,^*(\mathcal{L}_Q)$
is trivial.
\end{thm}

\begin{proof}
For the proof
it is enough to notice the basic facts relating with the Souriau map
and Fredholm pairs of Lagrangian subspaces(\cite{FO}, \cite{Le}).

For $\mu\in\mathcal{F}\Lambda_{\lambda}(H)$
let $p_{\lambda}(\mu)$ be 
\begin{equation*}
p_{\lambda}(\mu)
=\mathcal{P}(\mu^{\perp})+\mathcal{P}(\lambda^{\perp}),
\end{equation*}
then it is a positive Fredholm operator.
That is, we have a map
\begin{equation}
p_{\lambda}:\mathcal{F}\Lambda_{\lambda}(H)
\rightarrow \widehat{\mathcal{F}}_+(H).
\end{equation}
Then, 
$$
\Ker(p_{\lambda}(\mu))=\lambda\cap\mu,
$$
and
$$
\Coker(p_{\lambda}(\mu))=H/(\lambda^{\perp}+\mu^{\perp})
\cong \lambda/\left(\lambda\cap(\lambda^{\perp}+\mu^{\perp})\right).
$$

Also we know 
$$
\Ker (q_{\lambda}(\mu))=
\Ker(Id+\mathcal{S}_{\lambda}(\mu))=\lambda\cap\mu+J(\lambda\cap\mu)
\cong (\lambda\cap\mu)\otimes\mathbb{C}
$$
and since 
$Im ~(q_{\lambda}(\mu))
= \lambda\cap(\lambda\cap\mu)^{\perp}
+J(\lambda\cap(\lambda\cap\mu)^{\perp})$,
$$
\Coker(q_{\lambda}(\mu))= H_J/
\left(
\lambda\cap(\lambda\cap\mu)^{\perp}
+J(\lambda\cap(\lambda\cap\mu)^{\perp})\right)
\cong \left(\lambda/(\lambda\cap(\lambda^{\perp}+\mu^{\perp}))\right)\otimes\mathbb{C}.
$$

So the fiber of the induced bundle $q_{\lambda}^{\,*}(\mathcal{L}_Q)$ 
by the maps $q_{\lambda}$ is 
the complexification  
of that by the map $p_{\lambda}$, hence the bundle 
$q_{\lambda}\,^*(\mathcal{L}_Q)$ is trivial, since the Quillen
determinant is trivial on the subspace
$\widehat{\mathcal{F}}_+(H)
\subset\widehat{\mathcal{F}}_+(H\otimes\mathbb{C})$.
\end{proof}

\begin{cor}
The disjoint union
$$
\coprod
\limits_{\mu\in\mathcal{F}\Lambda_{\lambda}}
\stackrel{\dim \lambda\cap\mu} \wedge (\lambda\cap\mu)^*\otimes
\stackrel{\dim H/(\lambda^{\perp}+\mu^{\perp})}
\wedge H/(\lambda^{\perp}+\mu^{\perp})
$$
has a bundle structure as a (holomorphic)line bundle on the Fredholm Lagrangian Grassmannian
$\mathcal{F}\Lambda_{\lambda}(H)$ and is a trivial line bundle. 

Note that we do not have a particular
trivialization on the whole space $\mathcal{F}\Lambda_{\lambda}(H)$.
\end{cor}
\begin{rem}
For any Fredholm pair $(\lambda, \mu)$ of
Lagrangian subspaces
$$
\dim \lambda\cap\mu =\dim H/(\lambda^{\perp}+\mu^{\perp})
= \dim H/(\lambda+\mu).
$$
\end{rem}
\bigskip

Now let $\theta$ be a Lagrangian subspace which ``{\it almost coincides}'' with
$\lambda$ :
$$
\dim \lambda/(\lambda\cap\theta)= 
\dim \theta/(\lambda\cap\theta) <\infty.
$$
This relation is an equivalence relation among Lagrangian subspaces
and we denote it
by $\lambda\sim\theta$.  Then for such a pair $(\lambda,\theta)$,
$\lambda\sim\theta$, 
the Fredholm Lagrangian Grassmannian
coincides with each other :
$$
\mathcal{F}\Lambda_{\lambda}(H)
=\mathcal{F}\Lambda_{\theta}(H). 
$$

For a Lagrangian subspace $\theta$, let
us denote an open subset in $\mathcal{F}\Lambda_{\theta}(H)$
$$
\left
\{\mu\in\mathcal{F}\Lambda_{\theta}(H)~|\, ~\mu \cap\theta =\{0\}
\right\}
$$
by $\mathcal{F}\Lambda_{\theta}^{(0)}(H)$.
Then this space is isomorphic 
to the space of (real) bounded selfadjoint operators 
on $\theta$
and we have an open covering:
$$
\mathcal{F}\Lambda_{\lambda}(H)=\bigcup\limits_{\theta\sim\lambda}
\mathcal{F}\Lambda_{\theta}^{(0)}(H).
$$

On each open subset $\mathcal{F}\Lambda_{\theta}^{(0)}(H)$ 
($\theta\sim\lambda$)
we have a trivialization of the induced bundle
$q_{\lambda}\,^*(\mathcal{L}_Q)$
given by the trivialization (\ref{l-t-3})
on $\mathcal{U}_{A_{\theta}}$ with a trace class operator 
$A_{\theta}= -Id + \tau_{\theta}\circ\tau_{\lambda}$ (in fact this
is a finite rank operator).  Also there is a trivialization
on  $\mathcal{F}\Lambda_{\theta}^{(0)}(H)$ 
coming from the trivialization on an open subset 
$\mathcal{U}_{\mathcal{P}(\theta^{\perp})-\mathcal{P}(\lambda^{\perp})}\cap
\mathcal{F}(H)\subset\mathcal{F}(H\otimes\mathbb{C})$
through the map $p_{\lambda}$.
Here again the operator 
$\mathcal{P}(\theta^{\perp})-\mathcal{P}(\lambda^{\perp}) : H\to H$,
is a finite rank operator. For such two $\theta$ and
$\widetilde{\theta}$ ($\theta\sim\widetilde{\theta}$),
the transition function 
on the intersection $\mathcal{F}\Lambda_{\theta}^{(0)}(H)
\cap\mathcal{F}\Lambda_{\widetilde{\theta}}^{(0)}(H)$
is given
by the function through the map $q_{\lambda}$,
\begin{equation}
{\det}_{F}
\left\{(\tau_{\theta}-\tau_{\mu})
(\tau_{\widetilde{\theta}}-\tau_{\mu})^{-1}\right\}
={\det}_F\left\{(\mathcal{P}({\theta})-\mathcal{P}({\mu}))
(\mathcal{P}(\widetilde{\theta})-
\mathcal{P}({\mu}))^{-1}\right\},
\end{equation}
and that through the map $p_{\lambda}$ is
\begin{equation}
{\det}_F\left\{(\mathcal{P}(\theta^{\perp})+\mathcal{P}(\mu^{\perp}))
(\mathcal{P}(\widetilde{\theta}^{\perp})
+\mathcal{P}(\mu^{\perp}))^{-1}\right\}.
\end{equation}
Now we show these two functions coincide on
$\mathcal{F}\Lambda_{\theta}^{(0)}(H)
\cap\mathcal{F}\Lambda_{\widetilde{\theta}}^{(0)}(H)$ :
\begin{prop}
Let $\theta$ and $\widetilde{\theta}$ ``{\it almost coincide}'',
then for $\mu\in \mathcal{F}\Lambda_{\theta}^{(0)}(H)
\cap\mathcal{F}\Lambda_{\widetilde{\theta}}^{(0)}(H)$ we have 
\begin{align}
&{\det}_F\left\{(\mathcal{P}(\theta)-\mathcal{P}(\mu))
(\mathcal{P}(\widetilde{\theta})-
\mathcal{P}(\mu))^{-1}\right\}\label{1-term}\\
=&~
{\det}_F\left\{(\mathcal{P}(\theta^{\perp})-\mathcal{P}(\mu^{\perp}))
(\mathcal{P}(\widetilde{\theta}^{\perp})-
\mathcal{P}(\mu^{\perp}))^{-1}\right\}\label{eq-1}\\
=&
~{\det}_F\left\{(\mathcal{P}({\theta})+\mathcal{P}({\mu}))
(\mathcal{P}(\widetilde{\theta})
+\mathcal{P}(\mu))^{-1}\right\}\label{3-term}\\
=&
~{\det}_F\left\{(\mathcal{P}({\theta^{\perp}})+\mathcal{P}({\mu^{\perp}}))
(\mathcal{P}({\widetilde{\theta}^{\perp}})
+\mathcal{P}({\mu^{\perp}}))^{-1}\right\}.\notag
\end{align}
\end{prop}
\begin{proof}
Since $\mathcal{P}(x)= Id -\mathcal{P}(x^{\perp})$ for
any Lagrangian subspace $x$, 
we have
\begin{align*}
&(\mathcal{P}(\theta)-\mathcal{P}(\mu))
(\mathcal{P}(\widetilde{\theta})-
\mathcal{P}(\mu))^{-1}\\
=&~(\mathcal{P}(\theta^{\perp})-\mathcal{P}(\mu^{\perp}))
(\mathcal{P}(\widetilde{\theta}^{\perp})-
\mathcal{P}(\mu^{\perp}))^{-1}.
\end{align*}
This gives the first equality (\ref{eq-1}).

Next we prove the coincidence of the first term (\ref{1-term}) and
the third term (\ref{3-term}), then we know all the term coincide.

{}From the equality
\begin{align*}
&(\mathcal{P}({\theta})-\mathcal{P}({\mu}))
(\mathcal{P}({\widetilde{\theta}})-\mathcal{P}({\mu}))^{-1}\cdot
(\mathcal{P}({\widetilde{\theta}})-\mathcal{P}({\mu}))
(\mathcal{P}({\widetilde{\theta}})+\mathcal{P}({\mu}))^{-1}\\
=&
(\mathcal{P}({{\theta}})-\mathcal{P}({\mu}))
(\mathcal{P}({{\theta}})+\mathcal{P}({\mu}))^{-1}\cdot
(\mathcal{P}({\theta})+\mathcal{P}({\mu}))
(\mathcal{P}({\widetilde{\theta}})+\mathcal{P}({\mu}))^{-1}
\end{align*}
we have
\begin{align*}
&(\mathcal{P}({\theta})-\mathcal{P}({\mu}))
(\mathcal{P}({\widetilde{\theta}})-\mathcal{P}({\mu}))^{-1}\\
=&
(\mathcal{P}({{\theta}})-\mathcal{P}({\mu}))
(\mathcal{P}({{\theta}})+\mathcal{P}({\mu}))^{-1}\cdot
(\mathcal{P}({\widetilde{\theta}})-\mathcal{P}({\mu}))
(\mathcal{P}({\widetilde{\theta}})+\mathcal{P}({\mu}))^{-1}\\
&\cdot(\mathcal{P}({\widetilde{\theta}})-\mathcal{P}({\mu}))
(\mathcal{P}({\widetilde{\theta}})+\mathcal{P}({\mu}))^{-1}\cdot
(\mathcal{P}({\theta})+\mathcal{P}({\mu}))
(\mathcal{P}({\widetilde{\theta}})+\mathcal{P}({\mu}))^{-1}\\
&\cdot
(\mathcal{P}({\widetilde{\theta}})+\mathcal{P}({\mu}))
(\mathcal{P}({\widetilde{\theta}})-\mathcal{P}({\mu}))^{-1}.
\end{align*}
When we express the Lagrangian subspace $\theta\subset \mu+ J(\mu)= H$ 
as the graph of an operator $T_{\theta}:\mu\to J(\mu)$, the
operator
$(\mathcal{P}({\theta})-\mathcal{P}({\mu})) 
(\mathcal{P}({\theta})+\mathcal{P}({\mu}))^{-1}
$
is expressed in the following form:
$$
(\mathcal{P}({\theta})-\mathcal{P}({\mu})) 
(\mathcal{P}({\theta})+\mathcal{P}({\mu}))^{-1}
 :\begin{pmatrix}
x\\
J(y)\end{pmatrix}
\mapsto
\begin{pmatrix}
-Id&0\\
-T_{\theta}&Id
\end{pmatrix}
\begin{pmatrix}
x\\
J(y)
\end{pmatrix}, ~x,~y~\in~\mu.
$$
Hence we see that the operator
$$(\mathcal{P}({{\theta}})-\mathcal{P}({\mu}))
(\mathcal{P}({{\theta}})+\mathcal{P}({\mu}))^{-1}\cdot
(\mathcal{P}({\widetilde{\theta}})-\mathcal{P}({\mu}))
(\mathcal{P}({\widetilde{\theta}})+\mathcal{P}({\mu}))^{-1}
$$
is of the form:
$$
\begin{pmatrix}
Id&0\\
T_{\theta}-T_{\widetilde{\theta}}&Id
\end{pmatrix}.
$$

When $\theta$ and $\widetilde{\theta}$ almost coincide,
then this operator is of the form $Id$ $+$ ``{\it finite rank operator}'',
since $T_{\theta}-T_{\widetilde{\theta}}$ is a finite rank operator. 
Moreover we have 
$$
{\det}_F\begin{pmatrix}
Id & 0\\
T_{\theta}-T_{\widetilde{\theta}} & Id
\end{pmatrix}=1.
$$
Finally, together with an invariance of the Fredholm determinant
with respect to conjugations we have 
\begin{align*}
&{\det}_F\left\{(\mathcal{P}({\widetilde{\theta}})-\mathcal{P}({\mu}))
(\mathcal{P}({\widetilde{\theta}})+\mathcal{P}({\mu}))^{-1}\cdot
(\mathcal{P}({\theta})+\mathcal{P}({\mu})
(\mathcal{P}({\widetilde{\theta}})+\mathcal{P}({\mu}))^{-1}\right.\\
&\left.\qquad\qquad\qquad\qquad\cdot
(\mathcal{P}({\widetilde{\theta}})+\mathcal{P}({\mu}))
(\mathcal{P}({\widetilde{\theta}})-\mathcal{P}({\mu}))^{-1}\right\}\\
&\qquad\qquad={\det}_F\left\{(\mathcal{P}({\theta})+\mathcal{P}({\mu}))
(\mathcal{P}({\widetilde{\theta}})+\mathcal{P}({\mu}))^{-1}\right\},
\end{align*}
which proves the desired result.
\end{proof}

\begin{rem}
Although we know the triviality of the line bundle
$q_{\lambda}\,^*(\mathcal{L}_Q)$, there are no natural global
trivializations.  The {\it Maslov line bundle} 
on $\mathcal{F}\Lambda_{\lambda}$ (we do not define this here,
but is defined in a similar way as for the finite dimensional case,
see \cite{Ho})
is also a trivial line bundle just by its definition for which
the transition functions are given by 
the infinite dimensional analog of the 
{\it H\"ormander indexes}(\cite{FO}). 
So it is interesting
to give an isomorphism of these two line bundles on a particular subspace in 
the Fredholm Lagrangian Grassmannian in terms of a certain geometric
and/or analytic
data, which will give us a 
relation of the Fredholm determinant and
the Maslov index.
\end{rem}





\end{document}